\documentclass[12pt]{amsart}
\newcommand{\conservepaper}{
 \hoffset=-0.75in
 \setlength{\textwidth}{6.5in}
 \voffset=-0.5in
 \setlength{\textheight}{9.0in} 
 }
 \conservepaper

\usepackage{amssymb}

\newcommand{\integer}{{\mathbb Z}}
\newcommand{\real}{{\mathbb R}}
\newcommand{\rational}{{\mathbb Q}}
\renewcommand{\natural}{{\mathbb N}}
\newcommand{\G}{{\mathbb G}}

\newcommand{\F}{{\mathbb F}}

\newcommand{\comm}[2]{\lbrack\!\lbrack#1,#2\rbrack\!\rbrack}
\newcommand{\qn}{{\mathord{Q}}} 
\newcommand{\fs}{{\mathord{q}}} 
\newcommand{\pe}{{\mathord{r}}} 
\newcommand{\fdeg}{{\mathord{s}}} 
\newcommand{\Id}{\operatorname{Id}}
\newcommand{\field}{\mathord{F}}
\newcommand{\Fr}{\operatorname{Fr}}

\newcommand{\mdeg}{\mathop{\mathrm{deg}^-}}

\newcommand{\Gal}{\operatorname{Gal}}
\newcommand{\GL}{\operatorname{GL}}

\newcommand{\lead}[1]{\overline{#1}}

 \newcommand{\iso}{} 
\let\iso=\cong 
\renewcommand{\cong}{\equiv}

 \swapnumbers

\newcommand{\pref}[1]{{\upshape(\ref{#1})}}
\newcommand{\fullref}[2]{{\upshape\ref{#1}{\pref{#1-#2}}}}

\numberwithin{equation}{section}
\newtheorem{prop}[equation]{Proposition}
\newtheorem{thm}[equation]{Theorem}
\newtheorem{lem}[equation]{Lemma}
\newtheorem{cor}[equation]{Corollary}

\theoremstyle{definition}
\newtheorem{defn}[equation]{Definition}
\newtheorem{notation}[equation]{Notation}
\newtheorem{ack}[equation]{Acknowledgments}

\theoremstyle{remark}
\newtheorem{rem}[equation]{Remark}

 \newcounter{step}
\newenvironment{step}[1][\unskip]{\refstepcounter{step}
\em
 \medskip \noindent Step \thestep\
 #1.\ }{\unskip\upshape}

 \newcounter{case}
\newenvironment{case}[1][\unskip]{\refstepcounter{case}
\em
 \medskip \noindent Case \thecase\
 #1.\ }{\unskip\upshape}

\newcommand{\bigset}[2]{\left\{\, #1
 \mathrel{\left| \vphantom {\left\{ #1 \mid #2 \right\} } \right.}
 #2 \,\right\} }

\renewcommand{\see}[1]{(see~\ref{#1})}

\begin{document}

\title[Automorphisms of arithmetic subgroups]
 {On automorphisms of arithmetic subgroups of unipotent
groups in positive characteristic}

\author{Lucy Lifschitz}

\address{Department of Mathematics, University of
Virginia, Charlottesville, VA 22904}

\curraddr{Department of Mathematics, Tufts University, Medford, MA
02155}

\email{lifschitz@tufts.edu}

\author{Dave Witte}

\address{Department of Mathematics, Oklahoma State University,
Stillwater, OK 74078}

\email{dwitte@math.okstate.edu,
http://www.math.okstate.edu/$\sim$dwitte}

\thanks{Submitted to \emph{Transactions of the American Mathematical
Society} (July 2000). This version: 20 July 2000.} 

\begin{abstract}
 Let $\field$ be a local field of positive characteristic, and let $G$
be either a Heisenberg group over~$\field$, or a certain (nonabelian)
two-dimensional unipotent group over~$\field$. If $\Gamma$ is an
arithmetic subgroup of~$G$, we provide an explicit description of every
automorphism of~$\Gamma$. From this description, it follows that every
automorphism of~$\Gamma$ virtually extends to a virtual automorphism
of~$G$.
 \end{abstract}

\maketitle

\section{Introduction}

Roughly speaking, a discrete subgroup~$\Gamma$ of a topological
group~$G$ is automorphism rigid if every automorphism of~$\Gamma$
extends to a continuous automorphism of~$G$. However, the formal
definition below is slightly more complicated, because it allows for
passage to finite-index subgroups.

\begin{defn}
 It is traditional to say that a group~$\Gamma$ \emph{virtually} has a
property if some finite-index subgroup of~$\Gamma$ has the property.
It is convenient to extend this terminology to group isomorphisms.
 \begin{itemize}
 \item A \emph{virtual isomorphism} from~$G_1$ to~$G_2$ is an
isomorphism $\Lambda \colon G'_1 \to G'_2$, where $G'_i$ is a
finite-index, open subgroup of~$G_i$.
 \item A \emph{virtual automorphism} of~$G$ is a virtual isomorphism
from~$G$ to~$G$.
 \item A virtual isomorphism~$\Lambda$ from~$G_1$ to~$G_2$
\emph{virtually extends} an isomorphism $\lambda$
from~$\Gamma_1$ to~$\Gamma_2$ if there is a finite-index, open
subgroup~$\Gamma'_1$ of~$\Gamma_1$, such that $\Gamma'_1 \subset G_1$,
and $\Lambda|_{\Gamma'_1} = \lambda|_{\Gamma'_1}$.
 \end{itemize}
 \end{defn}

\begin{defn}
 A discrete subgroup~$\Gamma$ of a topological group~$G$ is
\emph{automorphism rigid} in~$G$ if every virtual automorphism
of~$\Gamma$ virtually extends to a virtual automorphism of~$G$.
 \end{defn}

A classical example is provided by the work of Malcev.

\begin{defn}[{\cite[Rem.~1.11, p.~21]{Raghunathan}}]
 A discrete subgroup~$\Gamma$ of a topological group~$G$ is a
(cocompact) \emph{lattice} if $G/\Gamma$ is compact.
 \end{defn}

\begin{thm}[{Malcev \cite{Malcev}, \cite[Cor.~2.11.1,
p.~34]{Raghunathan}}] \label{MalcevLie}
 If $\Gamma$ is a lattice in a $1$-connected, nilpotent real Lie
group~$G$, then $\Gamma$ is automorphism rigid in~$G$.
 
 In fact, every virtual automorphism of~$\Gamma$ extends to a unique
automorphism of~$G$.
 \end{thm}

Malcev's Theorem can be restated in the terminology of algebraic groups
(cf.\ \cite[after Thm.~2.12, p.~34]{Raghunathan}). Recall that a matrix
group~$G$ is \emph{unipotent} if, for every $g \in G$, there is some $n
\in \natural$, such that $(g - \Id)^n = 0$. (In other words, $1$~is the
only eigenvalue of~$g$.)

\begin{cor}
 Let $\Gamma$ be an arithmetic subgroup of a unipotent algebraic
$\rational$-group~$\G$. Then $\Gamma$ is an automorphism rigid lattice
in $\G(\real)$.
 \end{cor}

In this paper, we discuss the analogue of Malcev's Theorem for
unipotent groups over nonarchimedean local fields, instead of~$\real$.
It is well known that if $\G$ is a unipotent algebraic group over a
nonarchimedean local field~$L$ of characteristic zero, then the group
$\G(L)$ of $L$-points of~$\G$ has no nontrivial discrete subgroups.
(For example, $\integer$ is not discrete in the $p$-adic field
$\rational_p$.) Thus the case of characteristic zero is not of interest
in this setting; we will consider only local fields of positive
characteristic.

For abelian groups, it is easy to prove automorphism rigidity.

\begin{prop} \label{abelrigid}
 Let $\Gamma_1$ and~$\Gamma_2$ be lattices in a totally disconnected,
locally compact, abelian group~$G$. Then every isomorphism $\lambda
\colon \Gamma_1 \to \Gamma_2$ virtually extends to a virtual
automorphism~$\hat\lambda$ of~$G$.
 \end{prop}

\begin{proof}
 Since $\Gamma_1$ and~$\Gamma_2$ are
discrete, and $G$ is totally disconnected, there exists a compact, open
subgroup~$K$ of~$G$, such that $\Gamma_1 \cap K = \Gamma_2 \cap K = e$.
Let $\hat{G_1} = \Gamma_1 K$ and $\hat{G_2} = \Gamma_2 K$, so
$\hat{G_1}$ and~$\hat{G_2}$ are finite-index, open subgroups of~$G$,
and define $\hat\lambda \colon \hat{G_1} \rightarrow \hat{G_2}$ by
$\hat\lambda (\gamma c)=\lambda (\gamma ) \, c$ for $\gamma \in
\Gamma_1$ and $c \in K$. 
 \end{proof}

For nonabelian groups, automorphism rigidity seems to be surprisingly
more difficult to prove, but we provide examples of automorphism rigid
lattices. Although we do not have a general theory, and we do not have
enough evidence to support a specific conjecture, the examples suggest
that there may be mild conditions that imply arithmetic lattices are
automorphism rigid.

\begin{notation} \label{basicnotation}
 \begin{itemize}
 \item Fix a prime~$p$, and a power~$\fs$ of~$p$.
 \item $\F_\fs$ denotes the finite field of~$\fs$ elements.
 \item $\field$ denotes the field $\F_\fs((t))$ of formal power series
over~$\F_\fs$.
 \item $\field^-$ denotes $\F_\fs[t^{-1}]$, the $\F_\fs$-subalgebra
of~$\field$ generated by~$t^{-1}$.
 \end{itemize}
 Note that $\field$ is a local field of
characteristic~$p$. (Conversely, any local field of characteristic~$p$
is isomorphic to $\F_\fs((t))$, for some~$\fs$
\cite[Thm.~I.4.8, p.~20]{Weil}.) The subgroup~$\field^-$ is a lattice in
the additive group $(\field,+)$.
 \end{notation}

\begin{defn}
 Let $G$ be a closed subgroup of $\GL(m,\field)$, for some~$m \in
\natural$.
 \begin{itemize}
 \item Two discrete subgroups~$\Gamma_1$ and~$\Gamma_2$ of~$G$ are
\emph{commensurable} if $\Gamma_1 \cap \Gamma_2$ is a finite-index
subgroup of both $\Gamma_1$ and~$\Gamma_2$ \cite[p.~8]{MargBook}.
 \item A subgroup~$\Gamma$ of~$G$ is \emph{arithmetic} if it is
commensurable with $\GL(m,\field^-) \cap G$ (cf.\
\cite[\S I.3.1, pp.~60--62]{MargBook}).
 \end{itemize}
 By definition, if $\Gamma_1$ and~$\Gamma_2$ are arithmetic subgroups
of~$G$, then $\Gamma_1$ is commensurable with~$\Gamma_2$. Thus,
$\Gamma_1$ is a lattice in~$G$ if and only if $\Gamma_2$ is a lattice
in~$G$.
 \end{defn}

\begin{defn}[{cf.\ \cite[Ex.~9.2]{Borel-Springer-II}}]
\label{G2defn}
 Fix a power~$\pe$ of~$p$, and let
 $$ G_2 =
 \bigset{
 \begin{pmatrix}
 1 & y^\pe & z \\
 0 & 1   & y \\
 0 & 0   & 1
 \end{pmatrix}
 }{ y,z \in \field} .$$
 So $G_2$ is a two-dimensional, unipotent $\field$-group, and has
arithmetic lattices. Note that if $\pe > 1$, then $G_2$ is nonabelian.
 \end{defn}

The following theorem describes the virtual automorphisms of any
arithmetic lattice in~$G_2$. 

\begin{defn}
 For any continuous field automorphism~$\tau$ of~$\field$ and any $a \in \field
\setminus \{0\}$, there is a continuous automorphism~$\phi_{\tau, a}$
of~$G_2$, defined by
 $$ \phi_{\tau,a}
 \begin{pmatrix}
 1 & y^\pe & z \\
 0 & 1   & y \\
 0 & 0   & 1
 \end{pmatrix}
 =
 \begin{pmatrix}
 1 & a^\pe \tau(y)^\pe & a^{\pe+1} \tau(z) \\
 0 & 1   & a \tau(y) \\
 0 & 0   & 1
 \end{pmatrix}
 .$$
 Let us say that $\phi_{\tau,a}$ is \emph{standard} if 
 \begin{enumerate}
 \item there exist
$\sigma \in \Gal(\F_\fs/\F_p)$, $\alpha \in \F_\fs \setminus
\{0\}$, and $\beta \in \F_\fs$, such that
 $$ \tau \bigl( f(t^{-1}) \bigr) = \sigma \bigl( f( \alpha t^{-1} +
\beta) \bigr) ,$$
 for all $f(t^{-1}) \in \field$, and
 \item there exists some nonzero $b \in \field^-$, such that $a b \in
\field^-$.
 \end{enumerate}
 Note that if $\phi_{\tau,a}$ is standard, and $\Gamma$ is an
arithmetic lattice in~$G_2$, then $\phi_{\tau,a}(\Gamma)$ is
commensurable with~$\Gamma$.
 \end{defn}

\begin{thm} \label{autG2}
 Let 
 \begin{itemize}
 \item $\Gamma$ be an arithmetic lattice in~$G_2$; and
 \item $\lambda$ be a virtual automorphism of~$\Gamma$.
 \end{itemize}
 If $\pe > 2$, then there exist 
 \begin{itemize}
 \item a standard automorphism $\phi_{\tau,a}$ of~$G_2$, 
 \item a finite-index subgroup $\Gamma'$ of~$\Gamma$,
and 
 \item a homomorphism $\zeta \colon \Gamma' \to Z(\Gamma)$,
 \end{itemize}
 such that $\lambda(\gamma) = \phi_{\tau,a}(\gamma) \, \zeta(\gamma)$,
for all $\gamma \in \Gamma'$.
 \end{thm}

\begin{cor} \label{G2rigid}
 If $\pe \neq 2$, then any arithmetic lattice in~$G_2$ is automorphism
rigid.
 \end{cor}

Theorem~\ref{autG2} and Corollary~\ref{G2rigid} are proved
in Section~\ref{G2Section}. The authors do not know whether they remain true in the exceptional case
$\pe = p = 2$.

\begin{defn}
 Assume $p > 2$, let $\comm{\cdot}{\cdot} \colon \field^{2m} \times
\field^{2m} \to \field$ be a symplectic form, and, for notational
convenience, let $Z = F$. The corresponding \emph{Heisenberg group} is
the group
 $H = (\field^{2m} \times Z, \circ)$, where
 \[ ( v_1,z_1) \circ (v_2,z_2)= \bigl( v_1 + v_2, z_1 + z_2 +
 \comm{v_1}{v_2} \bigr) . \]
 We remark that, up to a change of basis, the symplectic form
$\comm{\cdot}{\cdot}$ on~$\field^{2m}$ is unique, so, up to isomorphism,
the Heisenberg group~$H$ is uniquely determined by~$m$. Note that $Z$ is
the center of~$H$.

Because $H$ is isomorphic to a subgroup of $\GL(m+2,\field)$, namely,
 $$ H \iso 
 \bigset{
 \begin{pmatrix}
 1 & x_1 & x_2 & \cdots & x_m & z \\
  & 1& & & & y_1 \\
  & & 1& \mbox{\Huge 0} & & y_2 \\
  & & & \ddots & & \vdots \\
  &\mbox{\Huge 0}& & & 1 & y_m \\
  & & & & & 1
 \end{pmatrix}
 }{ 
 \begin{matrix}
 x_1, \ldots, x_m \in \field, \\
 y_1, \ldots, y_m \in \field, \\
 z \in \field
 \end{matrix}
 } ,$$
 we may speak of arithmetic subgroups of~$H$. 
 
We assume that $\comm{\cdot}{\cdot}$ is defined over $\field^-$, by
which we mean that $\comm{\field^-}{\field^-} \subset \field^-$. Then
we may assume that the above isomorphism has been chosen so that 
 $$ \mbox{\emph{a subgroup~$\Gamma$ of~$H$ is arithmetic if and only if
it is commensurable with $(\field^-)^{2m} \times \field^-$.}} $$
 Thus, $H$ has arithmetic lattices.
 \end{defn}

We remark that one may define Heisenberg groups even if $p = 2$, but,
in this case, they are abelian, so they are not of particular interest.

\begin{defn}
 We say $T \in GL(2m,\field)$ is conformally symplectic if there exists
some nonzero $c_T \in \field$, such that, for all $v,w \in V$, we have
 $$\comm{T(v)}{T(w)} = c_T \, \comm{v}{w} .$$

 For every conformally symplectic $T \in GL(2m, \field)$, and every
 continuous field automorphism $\tau$ of $\field$, there is a
continuous automorphism $\phi _{T,\tau}$ of $H$ defined by 
$$ \phi_{T,\tau}(v,z) = \Bigl( \tau \bigl( T(v) \bigr), \tau( c_T z)
\Bigr) .$$
Let us say that $\phi_{T,\tau}$ is standard if
\begin{enumerate} \item there exist $\sigma \in \Gal(\F_q/\F_p)$,
$\alpha \in \F_q \setminus \{0\}$, and
  $\beta \in \F_q$, such that
 $$\tau \bigl( f(t^{-1}) \bigr) = \sigma \bigl( f(\alpha
  t^{-1} + \beta ) \bigr)$$ 
 for all $f(t^{-1}) \in \field$; and
 \item there exists some nonzero $b \in \field^-$, such that $b T \in
\operatorname{Mat}(2m, \field^-)$. \end{enumerate}
 \end{defn}

Note that if $\phi_{T,\tau}$ is standard, then $\phi_{T,\tau}(\Gamma)$
is commensurable with $\Gamma$ for any arithmetic lattice $\Gamma$
of~$H$.

\begin{thm}
\label{thm:heis}
 Assume $p > 2$. Let
 \begin{itemize}
 \item $\Gamma $ be an arithmetic lattice in a Heisenberg group~$H$; and
 \item $\lambda $ be a virtual automorphism of $\Gamma$.
 \end{itemize}
 Then there exist
 \begin{itemize}
 \item a standard automorphism $\phi _{T,\tau}$ of $H$;
 \item a finite index subgroup $\Gamma'$ of $\Gamma$; and
 \item a homomorphism $\zeta \colon \Gamma' \rightarrow Z(\Gamma)$,
 \end{itemize}
 such that $\lambda (\gamma) = \phi_{T,\tau}(\gamma) \, \zeta (\gamma
)$, for all $\gamma \in \Gamma'$.
\end{thm}  

\begin{cor}
\label{cor:heis}
 If $p > 2$, then any arithmetic lattice in a Heisenberg group~$H$ is
automorphism rigid.
 \end{cor}   

 Theorem~\ref{thm:heis} and
Corollary~\ref{cor:heis} are proved in Section~\ref{HeisSection}. 

\begin{rem}
 Malcev's Theorem~\ref{MalcevLie} does not extend to all lattices in
solvable Lie groups. (See the work of A.~Starkov \cite{Starkov} for a
thorough discussion.) On the other hand, the Mostow Rigidity Theorem
\cite{MostowRig} implies that lattices in most semisimple Lie groups
are automorphism rigid.

Superrigidity deals with extending homomorphisms, instead of only
isomorphisms. The Margulis Superrigidity Theorem \cite[Thm.~VII.5.9,
p.~230]{MargBook} implies that lattices in most semisimple Lie groups
are superrigid. (Lattices in many non-semisimple Lie groups are also
superrigid \cite{WitteSolvLatt}.) The Superrigidity Theorem also
applies to arithmetic subgroups of many semisimple groups defined over
nonarchimedean local fields, whether they are of characteristic zero or
not \cite{MargBook, Venkataramana}.
 \end{rem}

\begin{ack}
 The authors would like to thank the University of Bielefeld (Germany), 
the Isaac Newton Institute for Mathematical Sciences (Cambridge, U.K.),
the University of Virginia, and Oklahoma State University for their
hospitality. Most of this research was carried out during productive
visits to these institutions. 
 Financial support was provided by the German-Israeli Foundation for
Research and Development and the National Science Foundation
(DMS-9801136).
 \end{ack}

\section{Arithmetic subgroups of the two-dimensional unipotent
group~$G_2$} \label{G2Section}

Recall that $\pe$ and~$G_2$ are defined in Definition~\ref{G2defn}.
(Also recall the definitions of $p$, $\fs$, $\field$, and~$\field^-$ in
Notation~\ref{basicnotation}.)

\begin{proof}[{\bf Proof of Theorem~\ref{autG2}}]
 Let $\Gamma_1$ and~$\Gamma_2$ be finite-index subgroups of~$\Gamma$,
such that $\lambda$ is an isomorphism from~$\Gamma_1$ to~$\Gamma_2$.
 Then $\lambda$ induces isomorphisms 
 $$ \mbox{$\lambda^* \colon \Gamma_1/Z(\Gamma_1) \to
\Gamma_2/Z(\Gamma_2)$ and $\lambda_* \colon [\Gamma_1,\Gamma_1] \to
[\Gamma_2,\Gamma_2]$.} $$
 By identifying each of $G_2/Z(G_2)$ and $Z(G_2)$ with~$\field$ in the
natural way (and noting that $\Gamma_i \cap Z(G_2) = Z(\Gamma_i)$), we
may think of $\Gamma_i/Z(\Gamma_i)$ and $[\Gamma_i,\Gamma_i]$ as
$\F_p$-subspaces of~$\field$. By replacing $\Gamma_1$ and~$\Gamma_2$
with finite-index subgroups, we may assume that these subspaces are
contained in~$\field^-$. Then, because $\lambda$ is an isomorphism, we
see that the conditions of Notation~\ref{lambdas} are satisfied, so
Theorem~\ref{2D-lambda=at+b} below implies that there exist
 \begin{itemize}
 \item a standard automorphism~$\phi_{\tau,a}$ of~$G_2$, and
 \item a finite-index subgroup~$\Gamma_1'$ of~$\Gamma_1$,
 \end{itemize}
 such that $\lambda(\gamma) \in \phi_{\tau,a}(\gamma) \, Z(G)$, for
all $\gamma \in \Gamma_1'$. 

Because $\phi_{\tau,a}(\Gamma_1)$ is an arithmetic lattice, it is
commensurable with~$\Gamma_2$. Thus, replacing $\Gamma_1'$ with a
finite-index subgroup, we may assume that $\phi_{\tau,a}(\Gamma_1')
\subset \Gamma_2$. Then we may define $\zeta \colon \Gamma_1' \to
Z(\Gamma_2)$ by $\zeta(\gamma) = \lambda(\gamma) \,
\phi_{\tau,a}(\gamma)^{-1}$.
 \end{proof}

 \begin{lem} \label{lemma:ext}
 Let 
 \begin{itemize}
 \item $\Gamma$ be a lattice in a totally disconnected, locally compact
group~$G$,
 \item $A$ be a locally compact, abelian group, and 
 \item $\zeta \colon \Gamma \rightarrow A$ be a homomorphism. 
 \end{itemize}
 Assume 
 \begin{enumerate}
 \item there is a finite-index subgroup~$\Gamma'$ of~$\Gamma$, such
that $\Gamma ' \cap [G,G] \subset [\Gamma ,\Gamma ]$, and
 \item $\Gamma \cap [G,G]$ is a lattice in $[G,G]$.
 \end{enumerate}
 Then there is a finite-index, open subgroup~$\hat{G}$ of~$G$,
such that $\zeta$ extends to a continuous homomorphism
$\hat\zeta \colon \hat{G} \rightarrow A$ that is trivial on $[G,G]$.
 \end{lem}

\begin{proof}
 By assumption, there exists a lattice $\Gamma ' \subset \Gamma$ such
that $\Gamma ' \cap [G,G] \subset [\Gamma ,\Gamma ]$.  Since $\zeta
\colon \Gamma \rightarrow A$, and $A$ is abelian, we see that
$[\Gamma ,\Gamma ] \subset \ker \zeta $. Therefore $[\Gamma, \Gamma]
\subset \ker \zeta$, so, by the choice of~$\Gamma'$, we have
$\Gamma ' \cap [G,G] \subset \ker \zeta$.

By assumption, $\Gamma \cap [G,G]$ is a lattice in $[G,G]$, so $\Gamma
[G,G]/[G,G]$ is closed \cite[Thm.~1.13, p.~23]{Raghunathan}, hence
discrete. Thus, there is an open compact subgroup $K/[G,G] \subset
G/[G,G]$, such that $K \cap (\Gamma' [G,G])=e$. Let $\hat G = \Gamma' K
[G,G]$, and extend~$\zeta|_{\Gamma'}$ to a homomorphism $\hat\zeta
\colon \hat G' \rightarrow A$ by defining it to be trivial on $K
[G,G]$.
 \end{proof}

\begin{proof}[{\bf Proof of Corollary~\ref{G2rigid}}]
 We may assume $\pe > 2$. (Otherwise, we must have $\pe = 1$, which
means $G_2$ is abelian, so Proposition~\ref{abelrigid}
applies.)
 From Theorem~\ref{autG2}, we may assume there exist
 \begin{itemize}
 \item a standard automorphism $\phi_{\tau,a}$ of~$G_2$,  and 
 \item a homomorphism $\zeta \colon \Gamma_1 \to Z(\Gamma_2)$,
 \end{itemize}
 such that $\lambda(\gamma) = \phi_{\tau,a}(\gamma) \, \zeta(\gamma)$,
for all $\gamma \in \Gamma_1$. From
Lemma~\ref{lemma:ext}, we may assume that there is a finite-index
subgroup $G_2'$ of~$G_2$, such that $G_2'$ contains $[G_2,G_2]$,
and $\zeta$ extends to a homomorphism $\hat\zeta \colon G_2' \to
Z(G_2)$ that is trivial on $[G_2,G_2]$.  Let $G_2'' =
\phi_{\tau,a}(G_2')$.

Define $\hat\lambda \colon G_2' \to G_2$ by $\hat\lambda(g) =
\phi_{\tau,a}(g) \, \hat\zeta(g)$, for $g \in G_2'$, so $\hat\lambda$
is a continuous homomorphism that extends~$\lambda$. Because
$\hat\zeta$ is trivial on $[G_2,G_2]$, we know that
$\hat\lambda_{[G_2,G_2]} = \phi_{\tau,a}|_{[G_2,G_2]}$. Also, because
 $\hat\zeta(G_2') \subset Z(G_2) = [G_2,G_2]$, we know that
$\hat\lambda(g) \in \phi_{\tau,a}(g) \, [G_2,G_2]$ for all $g \in
G_2'$. Thus, $\hat\lambda$ induces an automorphism of $[G_2,G_2]$, and
an isomorphism $G_2'/[G_2,G_2] \to G_2''/[G_2,G_2]$, so $\hat\lambda$ is
an isomorphism.
 \end{proof}

\subsection{Using linear algebra to prove Theorem~\ref{autG2}}
 The remainder of this section is devoted to the statement and proof of
Theorem~\ref{2D-lambda=at+b}. This result is a reformulation of
Theorem~\ref{autG2} in terms of linear algebra. The reformulation is
not of intrinsic interest, but it clarifies the essential ideas of the
proof, and provides more flexibility, by allowing us to focus on the
important aspects of the internal structure of~$\Gamma$ that arise from
the structure of $\field^-$ as a polynomial algebra, without being
constrained by the external structure imposed by the group-theoretic
embedding of~$\Gamma$ in~$G_2$.

\begin{notation} 
 Define an $\F_p$-bilinear form $\comm{\cdot}{\cdot} \colon \field^- \times
\field^- \to \field^-$ by 
 $$\comm{a}{b} = a^\pe b - a b^\pe .$$
 For any $V,W \subset \field^-$, $\comm{V}{W}$ denotes the
$\F_p$-subspace of~$\field^-$ spanned by $\bigl\{ \, \comm{v}{w} \mid v
\in V, w \in W \, \bigr\}$.
 \end{notation}

 \begin{notation} \label{lambdas}
 Throughout the remainder of this section, we assume that
 \begin{itemize}
 \item $\pe > 2$;
 \item $V_1$ and $V_2$ are $\F_p$-subspaces of finite codimension
in~$\field^-$; and
 \item $\lambda^* \colon V_1 \to V_2$
and $\lambda_* \colon \comm{V_1}{V_1} \to \comm{V_2}{V_2}$ are
$\F_p$-linear bijections,
 \end{itemize}
 such that
 $$ \lambda_*\comm{a}{b} = \comm{\lambda^*(a)}{\lambda^*(b)} ,$$
 for all $a,b \in V_1$.
 \end{notation}

\begin{thm} \label{2D-lambda=at+b}
 There exist
 \begin{itemize}
 \item a subspace~$V_1'$ of finite codimension in~$V_1$,
 \item $a \in b^{-1} \field^-$, for some $b \in \field^-$,
 \item $\alpha,\beta \in \F_\fs$, with $\alpha \neq 0$, and 
 \item $\sigma \in \Gal(\F_\fs/\F_p)$, 
 \end{itemize}
 such that 
 $$ \lambda^* \bigl( f(t^{-1}) \bigr) = a \, \sigma \bigl( f( \alpha
t^{-1} + \beta) \bigr),$$
 for all $f(t^{-1}) \in V_1'$.
 \end{thm}

Let us outline the proof of Theorem~\ref{2D-lambda=at+b}, assuming, for
simplicity, that $V_1 = V_2 = \field^-$. For any power $\qn > 1$
of~$\pe$, we may define an equivalence relation on $\field^- \setminus
\{0\}$ by $a \mathbin{\equiv_\qn} b$ iff $a/b \in \field^\qn$; let
$[a]$ denote the equivalence class of~$a$. For each $a \in \field^-$,
the subspace $\comm{a}{\field^-}$ has infinite codimension in
$\comm{\field^-}{\field^-}$, but Proposition~\ref{<a,b>} shows that
$\comm{[a]}{\field^-}$ has finite codimension. Because
Corollary~\ref{lambda(a(k-)^q)} shows that
 $\lambda^*\bigl( [a] \bigr) = [ \lambda^*(a) ]$, 
 this codimension is a useful invariant.
Proposition~\ref{2D-codim[a,V]} shows that it is closely related to the
minimum degree of the elements of $[a]$. Using this,
Corollary~\ref{degpreserved} shows that there is some $a \in \field^-$,
a constant~$k$, and some~$\qn$, such that 
 $\mdeg \lambda^*(b) = k + \mdeg b$ for all $b \mathbin{\equiv_\qn} a$.
Also, Corollary~\ref{approxGCD} shows that $\lambda^*$ approximately
preserves the degrees of greatest common divisors. Then
Proposition~\ref{affineonq} shows that the restriction of~$\lambda^*$
to the $\F_p$-rational elements of some equivalence class is of the
desired form. Finally, we show that $\lambda^*$ has the desired form on
all of~$\field^-$.

\begin{notation}
 \begin{itemize}
 \item We use $\dim W$ to denote the dimension of a vector space~$W$
over~$\F_p$.
 \item Let $\fdeg = \dim \F_\fs$, so $\fs = p^\fdeg$.
 \item For $a = \sum_{i=0}^n \alpha_i t^{-i} \in \field^-$, with each
$\alpha_i \in \F_\fs$, we let $\mdeg a = n$ if $\alpha_n \neq 0$.
 \end{itemize}
 \end{notation}

The following proposition is used in almost all of the following
results. Because (\ref{<a,b>-[Gamma,Gamma]}~$\Rightarrow$) requires the
assumption that $e > 2$, it seems that a different approach will be
needed for the exceptional case $p = e = 2$.

\begin{prop} \label{<a,b>}
 \begin{enumerate}
 \item \label{<a,b>-[k,k]}
 The subspace $\comm{V_i}{V_i}$ has finite codimension in
$\field^-$.
 \item \label{<a,b>-[Gamma,Gamma]}
 Let $a,b \in V_i \setminus \{0\}$ and assume $a/b \notin \F_\fs$. The
subspace $\comm{a}{V_i} + \comm{b}{V_i}$ has finite codimension in
$\comm{V_i}{V_i}$ if and only if $a/b \in \field^\pe$.
 \end{enumerate}
 \end{prop}

\begin{proof}
 Because $\comm{a}{V_i}$ and $\comm{b}{V_i}$ have finite codimension in
$\comm{a}{\field^-}$ and $\comm{b}{\field^-}$, respectively, we see that
$\comm{a}{V_i} + \comm{b}{V_i}$ has finite codimension in
$\comm{a}{\field^-} + \comm{b}{\field^-}$. Thus, in
proving~\pref{<a,b>-[Gamma,Gamma]}, we may assume that $V_i =
\field^-$. 

\pref{<a,b>-[k,k]} This follows from our proof of
(\ref{<a,b>-[Gamma,Gamma]}~$\Leftarrow$) below.

 (\ref{<a,b>-[Gamma,Gamma]}~$\Leftarrow$) There are some nonzero $u,v
\in \field^-$, such that $a u^\pe = b v^\pe$. Let $x = a^\pe u - b^\pe
v$. 

We claim that $x \neq 0$. Otherwise, we
have 
 $$ a^{\pe^2 - 1} (a u^\pe) = (a^\pe u)^\pe = (b^\pe v)^\pe
 = b^{\pe^2 - 1} (b v^\pe) = b^{\pe^2 - 1} (a u^\pe) ,$$
 so $a^{\pe^2 - 1} = b^{\pe^2 - 1}$. This implies $a/b \in \F_\fs$,
which is a contradiction. This completes the proof of the claim. 

For any $y \in \field^-$, we have
 \begin{eqnarray*}
 \comm{a}{u y} - \comm{b}{v y} 
 &=& (a^\pe u y - a u^\pe y^\pe) - (b^\pe v y - b v^\pe y^\pe) \\
 &=& (a^\pe u y - b^\pe v y) - (a u^\pe y^\pe - b v^\pe y^\pe) \\
 &=& x y- 0 ,
 \end{eqnarray*}
 so $\comm{a}{\field^-} + \comm{b}{\field^-}$ contains $x \field^-$,
which is of finite codimension in~$\field^-$.

(\ref{<a,b>-[Gamma,Gamma]}~$\Rightarrow$) We may write $b$ (uniquely) in
the form
 $b = x + y^\pe a$, with $x,y \in \field$, and such that we
may write $x = \sum \alpha_i t^{-i}$ with $\alpha_i = 0$ whenever $i
\cong \mdeg(a) \pmod{\pe}$. (Note that we do \textbf{not} assume $x,y
\in \field^-$.)

For $u,v \in \field^-$, we have
 \begin{eqnarray*}
 \comm{a}{u} - \comm{b}{v}
 &=& (a^\pe u - a u^\pe) - (b^\pe v - b v^\pe) \\
 &=& (a^\pe u - b^\pe v) - \bigl( a u^\pe - (x + y^\pe a) v^\pe \bigr)
\\
 &=& (a^\pe u - b^\pe v) - a (u - y v)^\pe - x v^\pe .
 \end{eqnarray*}
 Whenever either $\mdeg(u)$ or $\mdeg(v)$ is large, it is obvious that
$\mdeg(a^\pe u - b^\pe v)$ is much smaller than 
 $ \max \bigl\{ \mdeg(u - y v)^\pe , \mdeg v^\pe \bigr\}$.
 Also, we may assume $x \neq 0$ (otherwise, we have $b/a = y^\pe \in
\field^\pe$, as desired), and, from the definition of~$x$, we know that
$\mdeg x \not\equiv \mdeg a \pmod{\pe}$, so
 $$ \mdeg \bigl( a (u - y v)^\pe - x v^\pe \bigr) 
 = \max \bigl\{ \mdeg \bigl( a (u - y v)^\pe\bigr) , \mdeg ( x v^\pe)
\bigr\} .$$
 Therefore, we conclude that 
 $$ \mdeg \bigl( \comm{a}{u} - \comm{b}{v} \bigr)
 \in \Bigl\{ \mdeg \bigl( a (u - y v)^\pe \bigr), \mdeg(x v^\pe)\Bigr\} 
 $$
 must be congruent to either $\mdeg(a)$ or $\mdeg(x)$, modulo~$\pe$.
Thus, because of our assumption that $\pe > 2$, we see that
$\comm{a}{\field^-} + \comm{b}{\field^-}$ does not contain elements of
all large degrees, so it does not have finite codimension
in~$\field^-$. Then, from~\pref{<a,b>-[k,k]}, we conclude that it does
not have finite codimension in $\comm{\field^-}{\field^-}$.
 \end{proof}

\begin{cor} \label{<a,c>&<b,c>}
 Let $a_1,a_2 \in V_i \setminus \{0\}$. We have $a_1/a_2 \in \field^\pe$
if and only if there is some nonzero $b \in V_1$, such that the subspace
$\comm{a_j}{V_i} + \comm{b}{V_i}$ has finite codimension in
$\comm{V_i}{V_i}$, for $j = 1,2$.
 \end{cor}

\begin{proof}
 ($\Rightarrow$) Choose $b \in a_1 \field^\pe \cap V_i \setminus
(\F_\fs a_1 \cup \F_\fs a_2)$. Then
Proposition~\fullref{<a,b>}{[Gamma,Gamma]} implies the desired
conclusion.

($\Leftarrow$) From Proposition~\fullref{<a,b>}{[Gamma,Gamma]}, we have
$a_1/b \in \field^\pe$ and $a_2/b \in \field^\pe$, so 
 $a_1/a_2 \in \field^\pe$.
 \end{proof}

\begin{lem} \label{k=ak^p}
 Let $a_1,a_2 \in \field^-$, and let $\qn > 1$ be a power of~$\pe$,
such that $\lambda^* \bigl( a_1 (\field^-)^\qn \cap V_1 \bigr) = a_2
(\field^-)^\qn \cap V_2$. Define
 \begin{itemize}
  \item subspaces $W_1$ and $W_2$ of finite codimension in~$\field^-$ by
 $a_i (\field^-)^\qn \cap V_i = a_i W_i^\qn$; 
 \item $\mu^* \colon W_1 \to W_2$ by 
 $\lambda^*(a_1 w^\qn) = a_2 \mu^*(w)^\qn$; and
 \item $\mu_* \colon \comm{W_1}{W_1} \to \comm{W_2}{W_2}$  by
 $\lambda_*(a_1^{\pe+1} w^\qn) = a_2^{\pe+1} \mu_*(w)^\qn$.
 \end{itemize}
 Then $\mu^*$ and~$\mu_*$ are $\F_p$-linear
bijections, and we have
 $$ \mu_*\comm{a}{b} = \comm{\mu^*(a)}{\mu^*(b)} ,$$
 for all $a,b \in W_1$.
 \end{lem}

\begin{defn}
 Let $\qn > 1$ be a power of~$p$.
 An element of~$\field^-$ is \emph{$\qn$-separable} if it is \textbf{not}
divisible by a nonconstant $\qn$th power.
 \end{defn}

\begin{cor} \label{lambda(a(k-)^q)}
 Let $a \in \field^-$, and let $\qn > 1$ be a power of~$\pe$, such that
$a$ is $\qn$-separable. Then there is some $\qn$-separable $b \in
\field^-$, such that $\lambda^* \bigl( a (\field^-)^\qn \cap V_1 \bigr)
= b (\field^-)^\qn \cap V_2$.
 \end{cor}

\begin{proof}
  Assume, for the moment, that $\qn = \pe$. For $a_1,a_2 \in \field^-
\setminus \{0\}$, define $a_1 \equiv a_2$ iff $a_1/a_2 \in \field^\pe$.
For nonzero $a,b \in V_1$, we see, from Notation~\ref{lambdas}, that
 $\comm{a}{V_1} + \comm{b}{V_1}$ has finite codimension in~$V_1$
 if and only if
 $\comm{\lambda^*(a)}{V_2} + \comm{\lambda^*(b)}{V_2}$ has finite
codimension in~$V_2$. Therefore, Corollary~\ref{<a,c>&<b,c>} implies
that $a \equiv b$ iff $\lambda^*(a) \equiv \lambda^*(b)$. The
equivalence classes are precisely the sets of the form $c(\field^-)^\pe
\cap V_i$, for some $\pe$-separable $c \in \field^-$, so the desired
conclusion is immediate.

We may now assume $\qn > \pe$. Let $\qn' = \qn/\pe$. There is some
$\qn'$-separable $a' \in \field^-$, such that $a \in a'
(\field^-)^{\qn'}$. By induction on~$\qn$, we know that there is some
$\qn'$-separable $b' \in \field^-$, such that
 $\lambda^* \bigl( a' (\field^-)^{\qn'} \cap V_1 \bigr) = b'
(\field^-)^{\qn'} \cap V_2$.

From the definition of~$a'$, we know there is some $a_1 \in \field^-$,
such that $a = a' a_1^{\qn'}$. Then, because $a$ is $\qn$-separable, we
know that $a_1$ is~$\pe$-separable. 

Define $W_1$, $W_2$, $\mu^*$, and~$\mu_*$ as in Lemma~\ref{k=ak^p}
(with $\qn'$, $a'$, and~$b'$ in the places of~$\qn$, $a$, and~$b$,
respectively). Because $a_1$ is $\pe$-separable, we know, from the case
$\qn=\pe$ in the first paragraph of this proof, that there is some
$\pe$-separable $b_1 \in \field^-$, such that 
 $\mu^* \bigl( a_1 (\field^-)^\pe \cap W_1 \bigr) = b_1 (\field^-)^\pe
\cap W_2$.
 Therefore
 \begin{eqnarray*}
 \lambda^* \bigl( a (\field^-)^\qn \cap V_1 \bigr)
 &=& \lambda^* \bigl[a' \bigl( a_1 (\field^-)^\pe \bigr)^{\qn'} \cap V_1
\bigr] \\
 &=& \lambda^* \bigl[ a' \bigl( a_1 (\field^-)^\pe \cap W_1
\bigr)^{\qn'} \bigr] \\
 &=& a' \bigl[ \mu^* \bigl( a_1 (\field^-)^\pe \cap W_1 \bigr)
\bigr]^{\qn'} \\
 &=& b' \bigl( b_1 (\field^-)^\pe \cap W_2 \bigr)^{\qn'} \\
 &=& b' \bigl( b_1 (\field^-)^\pe \bigr)^{\qn'} \cap V_2 \\
 &=& b' b_1^{\qn'} (\field^-)^\qn \cap V_2
 ,
 \end{eqnarray*}
 as desired.
  \end{proof}

\begin{lem} \label{2D-ideal}
 Let $a \in V_i$, let $\qn > 1$ be a power of~$\pe$, and let $k$ be the
codimension of~$V_i$ in~$\field^-$.
 Then there is some nonzero $b \in \field^-$ with $\mdeg b \le \pe^2
(k+1)$, such that $\comm{a (\field^-)^\qn \cap V_i}{V_i}$ contains a
codimension-$2k$ subspace of the ideal $a^\pe b^{\qn/\pe} \field^-$.
 \end{lem}

\begin{proof}
 Choose $c \in \field^- \setminus \F_\fs$, such that $a c^\qn \in V_i$
and $\mdeg c \le k+1$; let $b = c^{\pe^2} - c$.
  For $y \in \field^-$, we have
 \begin{eqnarray*}
 a^\pe b^{\qn/\pe} y
 &=& a^\pe (c^{\pe\qn} - c^{\qn/\pe}) y \\
 &=& (a^\pe c^{\pe\qn} y - a c^{\qn} y^\pe) - (a^\pe c^{\qn/\pe} y - a
c^{\qn} y^\pe) \\
 &=&  \comm{a c^{\qn}}{y} - \comm{a}{c^{\qn/\pe} y} \\
 &\in& \comm{a c^{\qn}}{\field^-} + \comm{a}{\field^-}
 ,
 \end{eqnarray*}
 so $\comm{a c^{\qn}}{\field^-} + \comm{a}{\field^-}$ contains $a^\pe
b^{\qn/\pe} \field^-$.

Because $\comm{a c^{\qn}}{V_i}$ and $\comm{a}{V_i}$ contain
codimension-$k$ subspaces of $\comm{a c^{\qn}}{\field^-}$ and
$\comm{a}{\field^-}$, respectively, this implies that $\comm{a
c^{\qn}}{V_i} + \comm{a}{V_i}$ contains a codimension-$2k$ subspace of
$a^\pe b^{\qn/\pe} \field^-$. Because both $a c^\qn$ and~$a$ belong to
$a (\field^-)^\qn \cap V_i$, the desired conclusion follows.
 \end{proof}

\begin{prop} \label{2D-codim[a,V]}
 Let $a \in V_i$, let $\qn > 1$~be a power of~$\pe$, and let $k$~be the
codimension of~$V_i$ in~$\field^-$. Then
 $$ \dim \frac{\field^-}{\comm{a(\field^-)^\qn \cap V_i}{V_i}}
 = \fdeg (\pe-1) (\mdeg a) + S + X ,$$
 where 
 \begin{itemize} 
 \item $S = \fdeg \max \{\, \mdeg c \mid \mbox{$c^\pe | a$, $c \in
\field^-$} \,\}$, and
 \item $0 \le X \le \fdeg \pe (k+1)\qn + 3k$.
 \end{itemize}
 \end{prop}

\begin{proof}
 Choose $b$ as in Lemma~\ref{2D-ideal}, and let $I = a^\pe b^{\qn/\pe}
\field^-$ and $\overline{\field^-} = \field^-/I$. It suffices to show
 \begin{equation} \label{2D-codim[a,V]-lowerbound}
 \dim \overline{ \field^- / \comm{a(\field^-)^\qn}{\field^-} } \ge \fdeg
(\pe-1) (\mdeg a) + S 
 \end{equation}
 and 
 \begin{equation} \label{2D-codim[a,V]-upperbound}
 \dim \overline{ \field^- / \comm{a}{\field^-} } \le S + \fdeg \pe^2
(k+1)\qn/\pe + \fdeg (\pe-1) \mdeg a .
 \end{equation}

 Let $u_1,u_2,\ldots,u_N$ be the irreducible factors of
 $a^\pe b^{\qn/\pe}$.
 Then we may write
 $$ a = u_1^{m_1} u_2^{m_2} \cdots u_f^{m_N} ,
 \qquad
 b^{\qn/\pe} = u_1^{\varepsilon_1} u_2^{\varepsilon_2} \cdots
u_f^{\varepsilon_N}  ,
 \mbox{ \qquad and \qquad\ }
 a^\pe b^{\qn/\pe} = u_1^{n_1} u_2^{n_2} \cdots u_f^{n_N} 
 ,$$
 where $n_j = \pe m_j+ \varepsilon_j$.

From the Chinese Remainder Theorem, we know that the natural ring
homomorphism from~$\overline{\field^-}$ to
 $$ \bigoplus_{j=1}^N \frac{\field^-}{ u_j^{n_j} \field^-} $$
 is an isomorphism.  Thus, we may work in each factor
$\field^-/u_j^{n_j} \field^-$, and add up the resulting codimensions.

 Define $\phi_j \colon \field^- \to \field^-/(u_j^{\pe m_j} \field^-)$
by
 $\phi_j(x) = a x^\pe$.
 Then, letting $m_j' = m_j - \lfloor m_j/\pe \rfloor$, we have
 $$ \ker \phi_j = \{\, x \in \field^- \mid u_j^{m_j'} | x \,\} ,$$
 so 
 \begin{eqnarray*}
 \dim \frac{\field^-}{u_j^{\pe m_j} \field^- + a (\field^-)^\pe} 
 &=& \dim \frac{\ker \phi_j}{u_j^{\pe m_j} \field^-} \\
 &=& \fdeg \dim_{\F_\fs} \frac{\ker \phi_j}{u_j^{\pe m_j} \field^-} \\
 &=& \fdeg (\pe m_j - m_j') \mdeg u_j \\
 &=& \fdeg (\pe - 1) \mdeg u_j^{m_j} + \fdeg \lfloor m_j/\pe\rfloor
\mdeg u_j .
 \end{eqnarray*}

 We have $a^\pe \in u_j^{\pe m_j} \field^-$, so
 \begin{equation} \label{2D-codim[a,V]-ak^Q}
 \comm{a (\field^-)^\qn}{\field^-}
 \subset a^\pe (\field^-)^{\qn \pe} \field^- + a (\field^-)^\qn
(\field^-)^\pe
 \subset u_j^{\pe m_j} \field^- + a (\field^-)^\pe
 \end{equation}
 and
 \begin{equation} \label{2D-codim[a,V]-a}
 \comm{a}{\field^-} + u_j^{\pe m_j} \field^-
 = u_j^{\pe m_j} \field^- + a (\field^-)^\pe .
  \end{equation}

 From \eqref{2D-codim[a,V]-ak^Q}, we have
 \begin{eqnarray*}
 \dim \frac{\field^-}{\comm{a (\field^-)^\qn}{\field^-} + u_j^{n_j}
\field^-}
 &\ge& \dim \frac{\field^-}{\comm{a (\field^-)^\qn}{\field^-} + u_j^{\pe
m_j} \field^-} \\
 &\ge& \dim \frac{\field^-}{u_j^{\pe m_j} \field^- + a (\field^-)^\pe}
\\
 &=& \fdeg (\pe - 1) \mdeg u_j^{m_j} + \fdeg \lfloor m_j/\pe\rfloor
\mdeg u_j,
 \end{eqnarray*}
 so
 \begin{eqnarray*}
 \dim \overline{ \field^- / \comm{a(\field^-)^\qn}{\field^-} } 
 &\ge& \sum_{j=1}^N \bigl( \fdeg (\pe-1) \mdeg u_j^{m_j} + \fdeg \lfloor
m_j/\pe\rfloor \mdeg u_j \bigr) \\
 &=& \fdeg (\pe-1) \mdeg a + S .
 \end{eqnarray*} 
 This establishes \eqref{2D-codim[a,V]-lowerbound}.

Because 
$\dim (u_j^{pm_j} \field^-/ u_j^{n_j} \field^-)
 = \fdeg \varepsilon_j \mdeg u_j$, and
from \eqref{2D-codim[a,V]-a},
 we have
 \begin{eqnarray*}
 \dim \frac{\field^-}{\comm{a}{\field^-} + u_j^{n_j} \field^-}
 &\le& \dim \frac{\field^-}{\comm{a}{\field^-} + u_j^{\pe m_j} \field^-}
+ \fdeg \varepsilon_j \mdeg u_j \\
 &=& \dim \frac{\field^-}{u_j^{\pe m_j} \field^- + a (\field^-)^\pe} +
\fdeg \varepsilon_j \mdeg u_j \\
 &=& \fdeg (\pe - 1) \mdeg u_j^{m_j} + \fdeg \lfloor m_j/\pe\rfloor
\mdeg u_j + \fdeg \varepsilon_j \mdeg u_j,
 \end{eqnarray*}
 so
 \begin{eqnarray*}
 \dim \overline{ \field^- / \comm{a}{\field^-} } 
 &\le& \sum_{j=1}^N \bigl(\fdeg (\pe - 1) \mdeg u_j^{m_j} + \fdeg
\lfloor m_j/\pe\rfloor \mdeg u_j + \fdeg \varepsilon_j \mdeg u_j\bigr)
\\
 &=& \fdeg (\pe-1) \mdeg a + S + \fdeg \mdeg b^{\qn/\pe} \\
 &\le&  \fdeg (\pe-1) \mdeg a + S + \fdeg \pe^2 (k+1) \qn/\pe
 .
 \end{eqnarray*}
 This establishes \eqref{2D-codim[a,V]-upperbound}.
 \end{proof}

\begin{lem} \label{a^p<a}
 For any $a \in \field^-$ and any $n \ge 0$, we have
 $$ \comm{a}{\field^-} + \comm{1}{\field^-} \subset
\comm{a^{\pe^n}}{\field^-} + \comm{1}{\field^-} .$$
 \end{lem}

\begin{proof}
 For any $v \in \field^-$, we have
 \begin{eqnarray*}
  \comm{v}{a} &=& v^\pe a - v a^\pe \\
 &=& v^\pe a - (v^{\pe^2} a^\pe - v^{\pe^2} a^\pe) - (v^\pe a^{\pe^2} -
v^\pe a^{\pe^2}) - v a^\pe \\
 &=& \comm{v^\pe a}{1} + \comm{v^\pe}{a^\pe} + \comm{v a^\pe}{1} \\
 &\in& \comm{1}{\field^-} + \comm{a^\pe}{\field^-} .
 \end{eqnarray*}
 Then the proof is completed by induction on~$n$.
 \end{proof}

\begin{prop} \label{deg(lambda(1))}
 There is some $N \in \natural$ {\upshape(}depending only on the
codimensions of $V_1$ and~$V_2$, not on the choice of $V_1$, $V_2$,
$\lambda^*$, or~$\lambda_*${\upshape)}, such that
$\mdeg \lambda^*(1)  \le N$.
 \end{prop}

\begin{proof}
 Let $k$ be the codimension of~$V_1$.
 Choose a power~$\qn > 1$ of~$\pe$ so large that $\lambda^*(1)$ is
$\qn$-separable. Then Corollary~\ref{lambda(a(k-)^q)} implies 
 $\lambda^* \bigl( (\field^-)^\qn \cap V_1 \bigr) = \lambda^*(1)
(\field^-)^\qn \cap V_2$.

 Choose $c \in \field^-\setminus \F_\fs$, such that $c^\qn \in V_1$ and
$\mdeg c \le r+1$.
 We have
 $$\begin{matrix}
 \comm{(\field^-)^{\qn} \cap V_1}{V_1}
 & \supset & \comm{1}{V_1} + \comm{c^{\qn}}{V_1} \hfill \\
 & \approx & \comm{1}{\field^-} + \comm{c^{\qn}}{\field^-} \hfill \\
 & \supset & \comm{1}{\field^-} + \comm{c^{\pe}}{\field^-} \hfill &
\mbox{\see{a^p<a}} \hfill \\
 & \supset & (c^{\pe^2} - c) \field^-  \hfill &
 \mbox{(proof of \pref{2D-ideal})} \hfill 
 .
 \end{matrix} $$
 So $\comm{(\field^-)^{\qn} \cap V_1}{V_1}$ has small codimension in
$\comm{V_1}{V_1}$.
 Therefore 
 $ \comm{\lambda^*(1) (\field^-)^{\qn} \cap V_2}{V_2} =
\lambda_*\comm{(\field^-)^{\qn} \cap V_1}{V_1}$
 must have small codimension in $\comm{V_2}{V_2}$, so
 $\mdeg \lambda^*(1)$ must be small, as desired.
 \end{proof}

\begin{cor} \label{almostseparable}
 There is some $N \in \natural$ {\upshape(}depending only on the
codimensions of $V_1$ and~$V_2$, not on the choice of $V_1$, $V_2$,
$\lambda^*$, or~$\lambda_*${\upshape)}, such that, for every power $\qn
> 1$ of~$\pe$ and every $\qn$-separable element~$a$ of~$V_1$, we have
$\mdeg \lambda^*(a) - \mdeg a' \le \qn N$, where $a'$ is the
$\qn$-separable element of $\lambda^*(a) \field^\qn$.
 \end{cor}

\begin{proof}
 Apply Proposition~\ref{deg(lambda(1))} to the map~$\mu^*$ of
Lemma~\ref{k=ak^p}.
 \end{proof}

\begin{prop} \label{qseparable}
 There is a power $\qn > 1$ of~$\pe$, and some $d > 0$, such that, for
every $v \in V_i$ with $\mdeg v > d$, there are $\qn$-separable
elements $v_1,\ldots,v_m$ of~$V_i$, such that
 $v = v_1 + \cdots + v_m$ and  
 $\mdeg v_j \le \mdeg v$, for $j = 1,\ldots,m$.
 \end{prop}

\begin{proof}
 Let $k$ be the codimension of~$V_i$ in~$\field^-$, and choose $\qn >
k+4$ so large that, for every $m \ge \qn$, the subspace~$V_i$ contains
elements of degree~$m$ whose leading coefficients span~$\F_\fs$. For
any element of~$V_1$ of degree~$m$, we show that there is a
$\qn$-separable element of~$V_i$ of degree~$m$ with the same leading
coefficient.

Let $\alpha$ be the leading coefficient of some element of~$V_1$ of
degree~$m$. Then $V_i$ contains exactly $\pe^{m-k}$ elements of
degree~$m$ with leading coefficient~$\alpha$.

On the other hand, if $a$ is an element of~$\field^-$ that
is of degree~$m$ and is not $\qn$-separable, then $a$ must be of the
form $a = x^\qn y$, where $x$ is an element of~$\field^-$ of some
degree~$j$, and $y$ is an element of~$\field^-$ of degree~$m-\qn j$.
Thus, the number of such elements~$a$ of degree~$m$ is no more than
 $$ \sum_{j=1}^\infty \fs^{j+1} \fs^{m-\qn j + 1}
 =  \fs^{m+2} \sum_{j=1}^\infty \fs^{j(1-\qn)}
 =  \frac{\fs^{m+2}}{\fs^{\qn-1}-1}
 \le  \frac{\fs^{m+2}}{\fs^{\qn-2}}
 <  \qn^{m-\qn+4}
 <  \frac{\qn^m}{\pe^k}
 .$$
 Therefore, not every element of~$V_i$ of degree~$m$ whose leading
coefficient is~$\alpha$ can be such an element~$a$, so $V_i$ has a
$\qn$-separable element of degree~$m$ with leading term~$\alpha$, as
desired.
 \end{proof}

\begin{cor} \label{deg(lambda(a))}
 For each $b \in \field^-$, there exists $N \in \natural$, such that,
for every $a \in b (\field^-)^\pe \cap V_1$, we have $|\mdeg
\lambda^*(a) - \mdeg a| \le N$.
 \end{cor}

\begin{proof}
 By symmetry, it suffices to show $\mdeg \lambda^*(a) \le \mdeg a + N$.
 We may assume $b$ is $\pe$-separable. By combining
Proposition~\ref{qseparable} with Lemma~\ref{k=ak^p}, we may choose a
power $\qn > 1$ of~$\pe$, such that each element of $b(\field^-)^\pe$
is a sum of $\qn$-separable elements of $b(\field^-)^\pe$ of smaller
degree. Thus, we may assume $a$ is $\qn$-separable (and our bound~$N$
may depend on~$\qn$).

Define $S$ as in the statement of Proposition~\ref{2D-codim[a,V]}, and
let $k_i$ be the codimension of~$V_i$. Because $a \in b (\field^-)^\pe$
and $b$ is $\pe$-separable, we have $S = \fdeg (\mdeg a - \mdeg
b)/\pe$, so Proposition~\ref{2D-codim[a,V]} implies
  $$ \left| \dim \frac{\field^-}{\comm{a(\field^-)^\qn \cap V_1}{V_1}}
 - \fdeg \bigl( \pe-1 + {\textstyle \frac{1}{\pe}} \bigr) \mdeg a
\right|
 \le \fdeg \frac{\mdeg b}{\pe} +  \bigl( \fdeg \pe (k_1+1)\qn + 3k_1
\bigr) $$
 is bounded.
 Similarly, letting $a'$ be the $\qn$-separable element of
$\lambda^*(a) \field^\qn$, and $b'$~be the $\pe$-separable element of
$\lambda^*(b) \field^\pe$, we know that
  $$ \left| \dim \frac{\field^-}{\comm{a'(\field^-)^\qn \cap V_2}{V_2}}
 - \fdeg \bigl( \pe-1 + {\textstyle \frac{1}{\pe}} \bigr) \mdeg a'
\right|
 \le \fdeg \frac{\mdeg b'}{\pe} +  \bigl( \fdeg \pe (k_2+1)\qn + 3k_2
\bigr) $$
 is bounded.
 Then, because 
 $$ \dim \frac{\comm{V_1}{V_1}}{\comm{a(\field^-)^\qn \cap V_1}{V_1}}
 = \dim \frac{\comm{V_2}{V_2}}{\comm{a'(\field^-)^\qn \cap V_2}{V_2}}
,$$
 we conclude that $|\mdeg a' - \mdeg a|$ is bounded.
Corollary~\ref{almostseparable} asserts that $|\mdeg \lambda^*(a) -
\mdeg a'|$ is also bounded.
 \end{proof}

\begin{cor} \label{degpreserved}
 For each $b \in \field^-$, there is a power~$\qn $ of~$\pe$, such that, for
every $a_1,a_2 \in b (\field^-)^\qn \cap V_1$, we have
 $\mdeg \lambda^*(a_1) - \mdeg \lambda^*(a_2) = \mdeg a_1 - \mdeg a_2$.
 \end{cor}

\begin{proof}
 Choose $N$ as in Corollary~\ref{deg(lambda(a))}. Now choose $\qn >
2N$. Because 
 $$\mdeg \lambda^*(a_1) \equiv \mdeg \lambda^*(a_2)  \pmod{\qn}
 \text{\qquad and\qquad} 
 \mdeg a_1 \equiv \mdeg a_2  \pmod{\qn},$$
 we have
 $$ \mdeg \lambda^*(a_1) - \mdeg a_1 \equiv \mdeg \lambda^*(a_2) - \mdeg
a_2 \pmod{\qn} ,$$
 so, from the choice of~$N$ and~$\qn$, we conclude that $\mdeg
\lambda^*(a_1) - \mdeg a_1 = \mdeg \lambda^*(a_2) - \mdeg a_2$.
 \end{proof}

\begin{prop}
 There is a constant $C > 0$, such that, for all $a_1,a_2 \in V_i$, and
every power~$\qn$ of~$\pe$, we have
 \begin{eqnarray*}
 \fdeg \mdeg \gcd(a_1,a_2) - C
 &\le& \dim \frac{\comm{V_i}{V_i}}{\comm{a_1(\field^-)^\qn \cap
V_i}{V_i} + \comm{a_2(\field^-)^\qn \cap V_i}{V_i}} \\
 &\le& C  \mdeg \gcd(a_1,a_2) + C
 .
 \end{eqnarray*}
 \end{prop}

\begin{proof}
 Because 
 $$\comm{a_1(\field^-)^\qn \cap V_i}{V_i} + \comm{a_2(\field^-)^\qn
\cap V_i}{V_i}
 \subset \gcd(a_1,a_2) \field^- ,$$
 the left-hand inequality is obvious.

Let $c = \gcd(a_1,a_2)$ and let $k$ be the codimension of~$V_i$. Then
Lemma~\ref{2D-ideal} implies that there exist nonzero $b_1,b_2 \in
\field^-$ with $\mdeg b_i \le \pe^2 (k+1)$, such that
 $\comm{a_j(\field^-)^\qn \cap V_i}{V_i}$
 contains a codimension-$2k$ subspace of $a_j^\pe b_j^{\qn/\pe}
\field^-$ for $j = 1,2$. Then, letting $b = b_1 b_2$, we have
$\mdeg b \le 2 \pe^2 (k+1)$, and
 $\comm{a_1(\field^-)^\qn \cap V_i}{V_i} + \comm{a_2(\field^-)^\qn \cap
V_i}{V_i}$
 contains a codimension-$4k$ subspace of the ideal $I = c^\pe
b^{\qn/\pe} \field^-$.
 Thus, it suffices to show that the codimension of
 $\comm{a_1}{\field^-} + \comm{a_2}{\field^-} + I$ in~$\field^-$ is
bounded above by $\fdeg (\pe+2) \mdeg c + \fdeg \mdeg b$.

Let $u_1,\ldots,u_N$ be the irreducible factors of 
 $c^\pe b^{\qn/\pe}$, so we may write
 $ c = u_1^{m_1} \cdots u_N^{m_N}$,
 $b = u_1^{\varepsilon_1} \cdots u_N^{\varepsilon_N}$, and
 $c^\pe b^{\qn/\pe} = u_1^{n_1} \cdots u_N^{n_N}$,
 where $n_j = \pe m_j + \varepsilon_j \qn/\pe$.
 From the Chinese Remainder Theorem, we have
 $\field^-/I \iso \bigoplus_{j=1}^N \field^-/ u^{n_j} \field^-$,
 so we may calculate the codimension in each factor, and then add them
up.

Fix~$j$. By interchanging $a_1$~and~$a_2$ if necessary, we may assume
that $u_j^{m_j+1} \nmid a_1$. It suffices to show that
 $$ \dim  \frac{\field^-}{\comm{a_1}{\field^-} + u^{n_j} \field^-}
 \le \fdeg \bigl( (\pe+2)m_j + \varepsilon_j \bigr) \mdeg u_j ;$$
 thus (because $m_j + \varepsilon_j \ge 1$), we need only show that
 $u_j^{(\pe+1)m_j+1} \field^- \subset \comm{a_1}{\field^-} + u_j^{n_j}
\field^-$.
 To show this, let $M$ be minimal, such that
 $u_j^{M+1} \field^- \subset \comm{a}{\field^-} + u_j^{n_j} \field^-$.
 (Obviously, we have $M < n_j$; we wish to show $M \le (\pe+1)m_j$.)
 Suppose $M > (\pe+1)m_j$. (This will lead to a contradiction.) We have
$m_j + \pe(M - \pe m_j) > M$, so
 \begin{eqnarray*} u_j^M \field^-
 &=& u_j^{\pe m_i} u_j^{M - \pe m_j} \field^- \\
 &\subset& a_1^\pe u_j^{M - \pe m_j} \field^- + u^{n_j} \field^- \\
 &\subset& \comm{a_1}{u_j^{M - \pe m_j} \field^-} + a_1 u_j^{\pe(M - \pe
m_j)} \field^- + u^{n_j} \field^- \\
 &\subset& \comm{a_1}{\field^-} + u_j^{m_j + \pe(M - \pe m_j)} \field^-
+ u^{n_j} \field^- \\
 &\subset& \comm{a_1}{\field^-} + u_j^{M+1} \field^- + u^{n_j} \field^-
\\
 &=& \comm{a_1}{\field^-} + u^{n_j} \field^-
 .
 \end{eqnarray*}
 This contradicts the minimality of~$M$.
 \end{proof}

\begin{cor} \label{approxGCD}
 There is a constant $C > 0$, such that, for all $a,b \in V_1$, we have
 $$ \frac{\mdeg \gcd(a,b)}{C} - C
 \le \mdeg \gcd \bigl( \lambda^*(a), \lambda^*(b) \bigr)
 \le C \mdeg \gcd(a,b) + C
 .$$
 \end{cor}

\begin{prop} \label{affineonq}
 There exist $b \in V_1$, $b' \in V_2$, $\alpha,\beta \in \F_\fs$, and
some~$\qn$ that is a power of both~$\pe$ and~$\fs$, such that, for all
$b \, f(t^{-\qn}) \in b \bigl( \F_p[t^{-1}] \bigr)^\qn \cap V_1$, we
have
 $\lambda^* \bigl( b \, f(t^{-\qn}) \bigr) = b' \, f( \alpha
t^{-\qn} + \beta)$.
 \end{prop}

\begin{proof}
 Corollary~\ref{degpreserved} shows that, by replacing~$V_1$ with some
$(\field^-)^\qn \cap V_1$ (using Lemma~\ref{k=ak^p}), we may assume
$\mdeg \lambda^*(a) = \mdeg a$, for every $a \in V_1$.

The terms $-C$ and~$+C$ in Corollary~\ref{approxGCD} are significant
only when $\mdeg \gcd(a,b)$ is small. On the other hand, $\mdeg
\gcd(a,b)$ can never be small (and nonzero) if $a,b \in (\field^-)^\qn$
for some large~$\qn$. Thus, by replacing $V_1$ with some $(\field^-)^\qn
\cap V_1$ (using Lemma~\ref{k=ak^p}), we may assume 
 $$ \frac{1}{C} \mdeg \gcd(a,b) \bigr)
 \le \mdeg \gcd \bigl( \lambda^*(a), \lambda^*(b) \bigr)
 \le C \mdeg \gcd(a,b) ,$$
 for every $a,b \in V_1$. In particular, $\gcd(a,b) = 1$ if and only if
 $\gcd \bigl( \lambda^*(a), \lambda^*(b) \bigr) = 1$.

Let $k$ be the codimension of~$V_1$ in~$\field^-$.
 Choose some 
 $$N > 4 \bigl( C (C+k) p^{C+k+1} + k + 1\bigr) .$$
 Choose a power~$\qn$ of~$\pe$, such that $\qn > N k$. There is some
nonzero $b \in \F_p[t^{-1}]$, with $\mdeg b \le N k$, such that
 $$ b (\F_p + t^{-\qn} \F_p + t^{-2\qn} \F_p + \cdots + t^{-N\qn} \F_p)
\subset V_1 .$$
 Because $\mdeg b < \qn$, we know that $b$ is $\qn$-separable, so, by
applying Lemma~\ref{k=ak^p} to $b (\field^-)^{\qn}\cap V_1$, we may
assume 
 $$ \F_p + t^{-1} \F_p + t^{-2} \F_p + \cdots + t^{-N} \F_p \subset V_1
.$$

By composing $\lambda^*$ with a map of the form $f(t^{-1}) \mapsto
\gamma f(\alpha t^{-1} + \beta)$, for some $\alpha, \beta, \gamma \in
\F_\fs$ (with $\alpha \gamma \neq 0$), we may assume $\lambda^*(1) = 1$
and $\lambda^*(t^{-1}) = t^{-1}$, so $\lambda^*|_{\F_p + \F_p t^{-1}} =
\Id$. 

Let $V_1^{\F_p} = V_1 \cap \F_p[t^{-1}]$. It suffices to show
$\lambda^*(a) = a$ for every $a \in V_1^{\F_p}$.

Suppose $\lambda^*|_{V_1^{\F_p}} \neq \Id$, and let 
 $$m = \min \bigset{ \mdeg a }{ \mbox{$\lambda^*(a) \neq a$, $a \in
V_1^{\F_p}$} } \ge 2 .$$
 Let $\Delta = \lambda^*(a) - a$, for any monic $a \in V_1^{\F_p}$ with
$\mdeg a = m$. (Note that the definition of~$m$ implies that $\Delta$ is
independent of the choice of~$a$.)

\setcounter{case}{0}

\begin{case}
 Assume $m \le N$.
 \end{case}
 Let $u$ be any irreducible element of~$\F_p[t^{-1}]$ with $\mdeg u \le
m-1$.

We claim that $V_1^{\F_p}$ contains a (monic) element~$a$, such
that $\mdeg a = m$ and $u | a$. To see this, let $b \in V_1^{\F_p}$ with $\mdeg b = m$. There is some $a \in \field^-$, such that
$u | a$ and $\mdeg(a-b) < \mdeg u < \mdeg b$. Because $\mdeg b \le N$,
this implies $a-b \in V_1^{\F_p}$, so $a \in V_1^{\F_p}$.

Because $u|a$ (and $\lambda^*(u) = u$), we know $\gcd \bigl( u,
\lambda^*(a) \bigr) \neq 1$. Because $u$ is irreducible, we conclude
that $u | \lambda^*(a)$. We also have $u | a$, so this implies $u |
(\lambda^*(a) - a) = \Delta$.

Thus, we see that $\Delta$ is divisible by every irreducible polynomial
over~$\F_p$ of degree $\le m-1$, so $\Delta$ is divisible by
$t^{-p^{m-1}} - t^{-1}$. Therefore $\mdeg \Delta \ge p^{m-1}$. However,
we also know $\mdeg \Delta \le \mdeg a = m$ (and all nonzero
polynomials in $\F_2[t^{-1}]$ are monic, so $\mdeg \Delta < m$ if $p =
2$). This is a contradiction.

\begin{case}
 Assume $m > N$.
 \end{case}
 Choose some monic $a \in V_1^{\F_p}$, with $\mdeg a = m$. By
subtracting a polynomial of degree $\le k$, we may assume $t^{-(k+1)} |
a$; let $u = a/ t^{-(k+1)}$. There is some nonzero $x \in \F_p[t^{-1}]$
with $\mdeg x \le k$, such that $u x \in V_1^{\F_p}$. (Note that $\mdeg
u x \le k + \mdeg u < m$.)

Let
 $$ \mathcal{C} = \{\, c \in \F_p[t^{-1}] \setminus \{0\} \mid \mdeg c <
C \,\} ,$$
 and
 $$b = \prod_{\mdeg c \le C+k} c, $$
 so $\mdeg b < (C+k) p^{C+k+1}$. 
 Now, for each $c \in \mathcal{C}$, let 
 $$ \mbox{$u_c = (u + c)x$ \qquad and \qquad $\displaystyle u'_c
= \frac{u_c}{\gcd(u_c,b)}$.}$$
 For $c \in \mathcal{C}$, we have $\{ c x, c t^{-(k+1)}\} \subset
V_1^{\F_p}$, so $u_c \in V_1^{\F_p}$ and $a + c t^{-(r+1)} \in
V_1^{\F_p}$. Also, because $a = u t^{-(k+1)}$, we have $(u + c) | (a +
c t^{-(k+1)})$. Then, since $\lambda^*(u
+ c) = u + c$, we have 
 $\mdeg \gcd \bigl( \lambda^*(a + c t^{-(k+1)}), u+c \bigr) \ge \bigl(
\mdeg (u+c) \bigr)/C $, so
  \begin{eqnarray*}
 \mdeg \gcd( \Delta, u'_c)
 &\ge&  \mdeg \gcd( \Delta, u_c) - \mdeg b \\
 &=&  \mdeg \gcd \bigl( \lambda^*(a + c t^{-(k+1)}) - (a + c
t^{-(k+1)}), u_c \bigr) - \mdeg b \\
 &\ge& \frac{\mdeg (u + c)}{C}  - \mdeg b \\
 &\ge& \frac{m-k-1}{C} - (C+k) p^{C+k+1} \\
 &\ge& \frac{m}{4C}
 .
 \end{eqnarray*}
 Also, for $c_1,c_2 \in \mathcal{C}$, we have 
 $$\mdeg \gcd(u_{c_1}, u_{c_2})
 \le \mdeg ( u_{c_1} - u_{c_2} )
 = \mdeg \bigl( (c_1-c_2)x \bigr)
 \le C + k ,$$
 so we see that $\gcd(u'_{c_1}, u'_{c_2}) = 1$ whenever
$c_1 \neq c_2$. Thus, we
conclude that
 $$ \mdeg \Delta \ge p^C \frac{m}{4C} > m .$$
 This is a contradiction.
 \end{proof}

\begin{proof}[{\bf Proof of Theorem~\ref{2D-lambda=at+b}}]
 Choose $b,b', \alpha, \beta,\qn$ as in Proposition~\ref{affineonq}. 
 By replacing $\lambda^*$ with $x \mapsto (b')^{-1} \,\lambda^*(bx)$
 and replacing $\lambda_*$ with $x \mapsto (b')^{-(\pe+1)}
\,\lambda^*(b^{\pe+1} x)$, we may assume $b = b' = 1$. Then, by
composing $\lambda^*$ and~$\lambda_*$ with $t^{-1} \mapsto
\alpha^{-1}(t^{-1} - \beta)$, we may assume $\alpha = 1$ and $\beta =
0$. Thus,
 \begin{equation} \label{lambda(a^Q)}
 \mbox{$\lambda^*(a) = a$
for all $a \in \F_p[t^{-\qn}] \cap V_1$.}
 \end{equation}
 We wish to show that there is some $\sigma \in \Gal(\F_\fs/\F_p)$,
such that, for every $a \in V_1$, we have $\lambda^*(a) = \sigma(a)$.

\setcounter{step}{0}

\begin{step} \label{sigmaforeach}
 For each $a \in V_1$, there is some $\sigma \in \Gal(\F_\fs/\F_p)$,
such that $\lambda^*(a) = \sigma(a)$.
 \end{step}
 Fix $a \in V_1$. Choose $C$ as in Corollary~\ref{approxGCD}, let $k$
be the codimension of~$V_1$, and choose $b \in \F_p[t^{-\qn}] \cap V_1$,
such that
 $$ \frac{\mdeg b}{C} - C > \qn \Bigl( \fdeg \bigl(\mdeg a + \mdeg
\lambda^*(a) \bigr) + k\Bigr) .$$
 Let 
 $$ c = \prod_{\sigma \in \Gal(\F_\fs/\F_p)} \bigl( b - \sigma(a) \bigr)
 \qquad \in \F_p[t^{-1}] ,$$
 and choose some nonzero $x \in \F_p[t^{-1}]$, such that $(c x)^\qn \in
V_1$ and $\mdeg x \le k$.

We have
 $$
 \begin{matrix}
 \qn \Bigl( \mdeg \gcd \bigl( b-\lambda^*(a) , c \bigr) + k \Bigr)
 &\ge& \mdeg \gcd \bigl( b-\lambda^*(a) , (c x)^\qn \bigr) \hfill \\
 &=& \mdeg \gcd \Bigl( \lambda^*(b - a) , \lambda^* \bigl( (c x)^\qn
\bigr) \Bigr) \hfill 
 & \mbox{\see{lambda(a^Q)}} \hfill \\
 &\ge& \displaystyle \frac{\mdeg \gcd \bigl( b - a ,  (c x)^\qn
\bigr)}{C} - C \hfill 
 & \mbox{(choice of~$C$)} \hfill \\
 &=& \displaystyle \frac{(\mdeg b)}{C} -C \hfill 
 & \mbox{($(b - a) | c$)} \hfill \\
 &>& \qn \Bigl( \fdeg \bigl( \mdeg a + \mdeg \lambda^*(a) \bigr) + k
\Bigr) \hfill 
 & \mbox{(choice of~$b$)} . \hfill 
 \end{matrix}
 $$
 Thus, from the definition of~$c$, we conclude that there is some
$\sigma \in \Gal(\F_\fs/\F_p)$, such that
 \begin{eqnarray*}
  \mdeg \gcd \bigl( b-\lambda^*(a) , b - \sigma(a)
\bigr)
 &>& \mdeg a + \mdeg \lambda^*(a) \\
 &=& \mdeg \sigma(a) + \mdeg \lambda^*(a) \\
 &\ge& \mdeg \bigl( \sigma(a) - \lambda^*(a) \bigr) \\
 &=& \mdeg \bigl( ( b-\lambda^*(a) ) - ( b - \sigma(a)) \bigr)
 .
 \end{eqnarray*}
 Therefore
 $\bigl ( b-\lambda^*(a) \bigr) - \bigl( b - \sigma(a) \bigr) = 0$, so
$\lambda^*(a) = \sigma(a)$.

\begin{step} \label{sigmaforFq}
 There is some $\sigma \in \Gal(\F_\fs/\F_p)$, such that $\lambda^*(a)
= \sigma(a)$ for every $a \in V_1$.
 \end{step}
 For $v \in \field^-$, let $\lead{v}$ denote the leading coefficient
of~$v$.
 Choose $b \in V_1$, such that $\lead{b}$ generates $\F_\fs$, that is,
$\F_\fs = \F_p[\lead{b}]$. From Step~\ref{sigmaforeach}, we know there
is some $\sigma \in \Gal(\F_\fs/\F_p)$, such that $\lambda^*(b) =
\sigma(b)$. We show $\lambda^*(a) = \sigma(a)$ for every $a \in V_1$.

Given $a \in V_1$, choose some $c \in V_1$, such that $\lead{c}$
generates $\F_\fs$, and such that $\mdeg c > \max\{ \mdeg a, \mdeg
b\}$. From Step~\ref{sigmaforeach}, there exist $\sigma', \sigma'' \in
\Gal(\F_\fs/\F_p)$, such that $\lambda^*(c) = \sigma'(c)$ and 
$\lambda^*(a+c) = \sigma''(a + c)$.
 Because $\mdeg c > \mdeg a$, we have $\lead{c} = \lead{a+c}$ and
$\lead{\lambda^*(a+c)} = \lead{\lambda^*(c)}$. Thus, we have
 $$ \sigma''(\lead{c})
 = \sigma''(\lead{a+c})
 = \lead{\sigma''(a+c)}
 = \lead { \lambda^*(a+c) }
 = \lead { \lambda^*(a) + \lambda^*(c) }
 = \lead { \lambda^*(c) }
 = \sigma'(\lead{c}) .$$
 Because $\lead{c}$ generates $\F_\fs$, we conclude that $\sigma'' =
\sigma'$. Therefore
 $$ \lambda^*(a) = \lambda^*(a+c) - \lambda^*(c)
 = \sigma''(a + c) - \sigma'(c)
 = \sigma'(a + c) - \sigma'(c)
 = \sigma'(a) .$$

 Similarly, we have $\lambda^*(b) = \sigma'(b)$. Because we also have
$\lambda^*(b) = \sigma(b)$, and $\lead{b}$ generates $\F_\fs$, we
conclude that $\sigma' = \sigma$. 

Therefore 
 $\lambda^*(a) = \sigma'(a) = \sigma(a)$, as desired.
 \end{proof}

\section{Arithmetic subgroups of Heisenberg groups} \label{HeisSection}

\begin{proof}[{\bf Proof of Theorem \ref{thm:heis}}]
Let $\Gamma _1$, $\Gamma _2$ be finite-index subgroups of $\Gamma $,
such that $\lambda \colon \Gamma _1 \rightarrow \Gamma _2$ is an
isomorphism. Let $\overline{\Gamma _i}$, $i=1,2$ be the image of $\Gamma _i$
in $\field ^{2m}$ under the projection $H \rightarrow \field ^{2m}$
with kernel $Z$. By passing to a finite-index subgroup, we can assume
that $\overline{\Gamma _i} \subset (\field^-)^{2m}$. Since $Z(\Gamma _i) =
\Gamma _i \cap Z$, we can identify $\overline{\Gamma _i}$ with $\Gamma
_i/Z(\Gamma _i)$, so $\lambda$ induces an isomorphism $\overline{\lambda
} \colon \overline{\Gamma _1} \rightarrow \overline{\Gamma _2}$.

\setcounter{step}{0}

\begin{step} \label{Heispf-linear}
 We can assume $\overline{\lambda }(av)=a \overline{\lambda }(v)$ for
all $a \in \field ^-$ and $v \in \overline{\Gamma _1}$, such that $av
\in \overline{\Gamma _1}$.
 \end{step}
 For each nonzero $v \in \overline{\Gamma _1}$, let $A_v=
\{\, a \in \field ^- \mid av \in \overline{\Gamma _1} \,\}$. Note that
$A_v$ is a finite-index subgroup of $\field ^-$. For $g,h \in
\Gamma_i$, we have $\field \overline{g} = \field \overline{h}$ if and
only if $C_{\Gamma_i}(g) = C_{\Gamma_i}(h)$, so
 $ \overline{\lambda }(A_v v) = \field \overline{\lambda }(v) \cap
\overline{\Gamma _2}$.
 Thus, we can define a function $\tau _v\colon A_v \rightarrow \field$
by $\tau _v(a) \overline{\lambda }(v)= \overline{\lambda }(av)$. Let $w \in
\overline{\Gamma _1}$ be such that $\comm{v}{w} \neq 0$, and let $a
\in A_v \cap A_w$. Then
 \begin{eqnarray*}
 \tau _v(a) \comm{\overline{\lambda }(v)}{\overline{\lambda }(w)}
 &=& \comm{\overline{\lambda }(av)}{\overline{\lambda }(w)} \\
 &=& \lambda \bigl( \comm{av}{w} \bigr) \\
 &=& \lambda \bigl( \comm{v}{aw} \bigr) \\
 &=& \comm{\overline{\lambda }(v)}{\overline{\lambda }(aw)} \\
 &=& \tau _w(a) \comm{\overline{\lambda }(v)}{\overline{\lambda }(w)}
 .
\end{eqnarray*}
Thus 
\begin{equation}
\label{eqn:eq}
 \mbox{$\tau _v = \tau _w$  on $A_v \cap A_w$ whenever $\comm{v}{w}
\neq 0$.}
 \end{equation}

For any nonzero $v,w \in \overline{\Gamma _1}$ and any $a \in A_v \cap A_w$,
since $\overline{\Gamma _1} \cap a^{-1}\overline{\Gamma _1}$ is of finite index in
$\overline{\Gamma _1}$, we can find $u \in \overline{\Gamma _1}$ so that $a \in
A_u$, $\comm{u}{v} \neq 0$, and $\comm{u}{w} \neq 0$. Then it follows
from Equation~\pref{eqn:eq} that $\tau _v(a) =\tau _u(a)=\tau _w(a)$.
Since $a \in A_v \cap A_w$ was arbitrary, we conclude that 
\begin{equation} \label{eqn:eqall}
 \mbox{$\tau _v =
\tau _w$ on $A_v \cap A_w$, for all nonzero $v,w \in \overline\Gamma_1$.} 
 \end{equation}
 For an arbitrary $a \in \field ^-$ we can always find $w \in
\overline{\Gamma_1}$ so that $a \in A_w$, thus we can define a function
$\tau \colon \field^- \rightarrow F$, by $\tau (a) = \tau _w(a)$.
Equation~\pref{eqn:eqall} implies that $\tau$ is well defined. Note
that $\tau (1)=1$. Since
 \begin{eqnarray*}
 \tau (a)\tau (b)\comm{\overline{\lambda }(u)}{\overline{\lambda }(v)} 
 &=& \comm{\overline{\lambda }(au)}{\overline{\lambda }(bv)} \\
 &=& \lambda \bigl( \comm{au}{bv} \bigr) \\
 &=& \lambda \bigl( \comm{abu}{v} \bigr) \\
 &=& \comm{\overline{\lambda }(abu)}{\overline{\lambda }(v)} \\
 &=& \tau (ab)\comm{\overline{\lambda }(u)}{\overline{\lambda }(v)} ,
 \end{eqnarray*}
 we have $\tau (a)\tau (b)=\tau (ab)$. Since $\tau $ is also an
additive homomorphism, and $\overline{\lambda }$ is an isomorphism, we
conclude that $\tau $ is a ring automorphism of~$\field ^-$. Therefore
$\tau \bigl( f(t^{-1}) \bigr) = \sigma \bigl( f(\alpha t^{-1} + \beta )
\bigr)$ for $f(t^{-1}) \in \field ^-$, where $\sigma \in
\Gal(\F_q/\F_p)$, $\alpha \in \F_q \setminus \{0\}$, and $\beta \in
\F_q$. Hence, by composing with the standard automorphism
$T_{\Id,\tau^{-1}}$, we obtain the claim.

\begin{step} \label{Heispf-Z}
 We may assume that $\lambda|_{Z(\Gamma _1)}$ is the identity map.
 \end{step}
 Let $v_1,w_1,v_2,w_2 \in \overline{\Gamma}$ with $\comm{v_i}{w_i} \neq 0$.
There is a finite-index subgroup~$A$ of~$\field^-$, such that $av_i \in
\overline{\Gamma}$, for every $a \in A$ and $i = 1,2$. Then, for all $a \in
A$, Step~\ref{Heispf-linear} implies that
 $$ \frac{\lambda \bigl( a \comm{v_i}{w_i} \bigr))}{ a \comm{v_i}{w_i}}
 =  \frac{\lambda \bigl( \comm{v_i}{w_i} \bigr)}{\comm{v_i}{w_i}} .
 $$
 Thus, choosing $a_1,a_2 \in A$, such that 
 $a_1 \comm{v_1}{w_1} = a_2 \comm{v_2}{w_2}$, we have
 $$ \frac{\lambda \bigl( \comm{v_1}{w_1} \bigr)}{\comm{v_1}{w_1}}
 = \frac{\lambda \bigl( a_1 \comm{v_1}{w_1} \bigr)}{a_1 \comm{v_1}{w_1}}
 = \frac{\lambda \bigl( a_2 \comm{v_2}{w_2} \bigr)}{a_2 \comm{v_2}{w_2}}
 = \frac{\lambda \bigl( \comm{v_2}{w_2} \bigr)}{\comm{v_2}{w_2}}
 .$$
 We conclude that $\lambda(z)/z = C$ is constant,
for $z \in \comm{\overline{\Gamma_1}}{\overline{\Gamma_1}} \setminus \{0\}$.

By composing with a standard automorphism~$\phi_{T,\Id}$, such that
$c_T = 1/C$, we may assume that $C = 1$, so
$\lambda|_{[\Gamma_1,\Gamma_1]} = \Id$. Then, by replacing $\Gamma_1$
with a finite-index subgroup~$\Gamma_1'$, such that $\Gamma_1' \cap Z
\subset [\Gamma_1,\Gamma_1]$, we may assume $\lambda|_{Z(\Gamma_1)} =
\Id$.

\begin{step} \label{Heispf-symp}
 $\overline{\lambda }\colon \overline{\Gamma _1} \rightarrow
\overline{\Gamma _1}$ can be extended to a conformally symplectic map
$\overline{\Lambda} \colon \field ^{2m} \rightarrow \field ^{2m}$, with
$c_{\overline{\Lambda }}=1$.
 \end{step}
 By Step~\ref{Heispf-linear}, $\overline{\lambda }(av)= a
\overline{\lambda }(v)$ for all $a \in \field ^-$ and $v \in
\overline{\Gamma _1}$ such that $av \in \overline{\Gamma_1}$. Because
$\overline{\Gamma_1}$ is commensurable with $(\field^-)^{2m}$, this
implies that $\overline{\lambda}$ extends (uniquely) to an
$\field$-linear map $\overline{\Lambda }\colon \field ^{2m} \rightarrow
\field ^{2m}$. For any $v,w \in \overline{\Gamma_1}$, we have
 $$ \comm{\overline{\Lambda }(v)}{\overline{\Lambda }(w)}
 = \comm{\overline{\lambda }(v)}{\overline{\lambda }(w)}
 = \lambda \bigl( \comm{v}{w} \bigr)
 =  \comm{v}{w} ,$$
 by Step~\ref{Heispf-Z}. Because $\overline{\Gamma_1}$
spans~$\field^{2m}$, this implies that $\overline{\Lambda}$ is
conformally symplectic, with $c_{\overline{\Lambda }}=1$. 

 \begin{step}
 Completion of the proof.
 \end{step}
 Define $\hat{\Lambda } \colon H \rightarrow H$ by $ \hat{\Lambda
}(v,z)= \bigl( \overline{\Lambda }(v),z \bigr)$.  From
Step~\ref{Heispf-symp}, we see that $\hat{\Lambda }$ is an automorphism.
 Denote by $\zeta \colon \Gamma_1 \rightarrow Z(H)$ the map defined
by $\zeta (\gamma )=\hat{\Lambda }(\gamma )^{-1} \lambda (\gamma )$.
Then $\zeta $ is a homomorphism and $\lambda (\gamma )= \zeta (\gamma )
\, \hat{\Lambda }(\gamma )$, for $\gamma \in \Gamma_1$.
 \end{proof}

\begin{proof}[{\bf Proof of Corollary \ref{cor:heis}}]
 From Theorem \ref{thm:heis}, we may assume there exist
 \begin{itemize}
 \item a standard automorphism $\phi _{T,\tau}$ of $H$; and
 \item a homomorphism $\zeta\colon \Gamma _1 \rightarrow Z(H)$,
 \end{itemize}
 such that $\lambda (\gamma )=\phi _{T,\tau}(\gamma ) \, \zeta
(\gamma)$ for all $\gamma \in \Gamma _1$.
 By Lemma~\ref{lemma:ext}, there exists a finite-index open subgroup
$\hat{H}$ of~$H$, containing $[H,H]$, such that $\zeta$ extends to
$\hat{\zeta } \colon \hat{H} \rightarrow Z(H)$. Let $H' =
\phi_{T,\tau}(\hat{H})$.

Define $\hat{\Lambda } \colon \hat{H} \rightarrow H$ by
$\hat{\Lambda}(h)= \phi _{T,\tau}(h) \, \hat{\zeta }(h)$, so that
$\hat{\Lambda }$ is a continuous homomorphism virtually extending
$\lambda $. Because $\hat{\zeta }$ is trivial on $[H,H]$, we have
$\hat{\Lambda }|_{[H,H]} = \phi _{T,\tau}|_{[H,H]}$, so
$\hat{\Lambda }|_{[H,H]}$ is an automorphism. Because
$\hat{\zeta}(\hat{H}) \subset Z(H) = [H,H]$, we see that
$\hat{\Lambda}$ induces an isomorphism $\hat{H}/[H,H] \rightarrow
H'/[H,H]$. So $\hat{\Lambda} \colon \hat{H} \rightarrow H'$ is an
isomorphism.  
 \end{proof}

\begin{defn}
 Let
 $$ H_p = \bigset{
 \begin{pmatrix}
 1 & x_1^p & x_2^p & \cdots & x_m^p & z \\
  & 1& & & & y_1^p \\
  & & 1& \mbox{\Huge 0} & & y_2^p \\
  & & & \ddots & & \vdots \\
  &\mbox{\Huge 0}& & & 1 & y_m^p \\
  & & & & & 1
 \end{pmatrix}
 }{
  \begin{matrix}
 x_1, \ldots, x_m \in \field, \\
 y_1, \ldots, y_m \in \field, \\
 z \in \field
 \end{matrix}
 }
 .$$ 
 \end{defn}

\begin{rem} $H_p$ could also be described as the $\field$-points of the
group obtained from $H$ by applying the isogeny of factoring by the Lie
algebra of $Z(H)$ \cite[Prop.~V.17.4, p.~215]{Borel}.
 \end{rem}

\begin{cor}
 Any arithmetic lattice in~$H_p$ is automorphism rigid.
 \end{cor} 

\begin{proof}
 Let $\lambda_p \colon \Gamma_1 \to \Gamma_2$ be an isomorphism, where
$\Gamma_1$ and~$\Gamma_2$ are arithmetic lattices in~$H_p$.
 Define 
 $$H'_p = 
 \bigset{
 \begin{pmatrix}
 1 & x_1^p & x_2^p & \cdots & x_m^p & z^p \\
  & 1& & & & y_1^p \\
  & & 1& \mbox{\Huge 0} & & y_2^p \\
  & & & \ddots & & \vdots \\
  &\mbox{\Huge 0}& & & 1 & y_m^p \\
  & & & & & 1
 \end{pmatrix}
 }{
  \begin{matrix}
 x_1, \ldots, x_m \in \field, \\
 y_1, \ldots, y_m \in \field, \\
 z \in \field
 \end{matrix}
 }
 $$
 and
 $$ A =
 \bigset{
 \begin{pmatrix}
 1 & 0 & 0 & \cdots & 0 & z \\
  & 1& & & & 0 \\
  & & 1& \mbox{\Huge 0} & & 0 \\
  & & & \ddots & & \vdots \\
  &\mbox{\Huge 0}& & & 1 & 0 \\
  & & & & & 1
 \end{pmatrix}
 }{
  \begin{matrix}
 z = \displaystyle
 \sum_{ \begin{matrix}
             0\le i \le n \\
             i \not\equiv 0 \pmod{p}
        \end{matrix}}
 \alpha_i t^{-i}, \\
 \\
 n \in \natural, \\
 \alpha_i \in \F_\fs \\
 \end{matrix}
 }
 $$
 Then $H_p=H'_p \times A$. By passing to a finite-index subgroup we may
assume that $\Gamma _1=\Gamma '_1 \times \Gamma _{1,A}$, where $\Gamma
'_1 = \Gamma _1 \cap H_p'$ and $\Gamma _{1,A}=\Gamma _1 \cap A$. Let
$\Omega = \lambda _p(\Gamma _{1,A}) \subset Z(\Gamma _2)$ and $\Gamma
'_2 = \lambda _p(\Gamma '_1)$. Then, by passing to a finite-index
subgroup, we may assume $\Omega \cap H'_p = e$ and $\Gamma '_2 \cap A
=e$.

\setcounter{step}{0}

\begin{step} 
 Let $\pi _A:Z(H_p) \rightarrow A$ denote the projection with kernel $H'_p$.
 Then $\pi _A \circ \lambda _p: \Gamma _{1,A} \rightarrow \pi
 _A(\Omega )$ virtually extends to a virtual automorphism $\Psi $ of $A$.
\end{step}
 It is easy to see that 
$\pi _A(Z(\Gamma _2))$ is closed in $A$ and hence is a lattice. Because
$Z(\Gamma '_1) \times \Gamma _{1,A}$ has finite index in $Z(\Gamma
_1)$, we know $\lambda _p(Z(\Gamma '_1)) \times \lambda _p(\Gamma
_{1,A})$ has finite index in $Z(\Gamma _2)$. Then, since $[\Gamma '_1,
\Gamma '_1]$ has finite index in $Z(\Gamma '_1)$ and
\[ \lambda _p([\Gamma '_1,\Gamma '_1]) \subset [\Gamma '_2,\Gamma '_2]
\subset H'_p = \ker \pi _A \]
we conclude that $\pi _A(\Omega )=\pi _A(\lambda _p(\Gamma _{1,A}))$
has finite index in $\pi _A(Z(\Gamma _2))$. Hence $\pi _A(\Omega )$ is
a lattice in $A$. By Proposition \ref{abelrigid} $\pi _A \circ
\lambda _p: \Gamma _{1,A} \rightarrow \pi _A(\Omega )$ virtually
extends to a virtual automorphism $\Psi $ of $A$.

 \begin{step}
  Let $\pi' \colon H_p \to H_p'$ be the projection with kernel~$A$, and
let
 $\mu_p = \pi' \circ \lambda_p|_{\Gamma_1'} \colon \Gamma_1' \to
\pi'(\Gamma_2')$.
 Then $\mu_p$ virtual extends to a virtual automorphism of~$H_p'$.
 \end{step}
 
We claim that $\pi'(\Gamma_2')$ is an arithmetic lattice in~$H_p'$.
 Because $\Gamma_1 = \Gamma_1' \times \Gamma_{1,A}$ and $\Gamma_{1,A}
\subset Z(\Gamma_1)$, we have
 $$ \Gamma_2 = \Gamma_2' \times \Omega \subset \Gamma_2' \, Z(H_p) .$$
 Then, because $\Gamma_2' \subset \Gamma_2$, we conclude that
 $\Gamma_2' \, Z(H_p) = \Gamma_2' \, Z(H_p)$
 is a lattice in $H_p/Z(H_p) \iso H_p'/ Z(H_p')$.
 So the image of $\pi'(\Gamma_2')$ in $H_p'/ Z(H_p')$ is a lattice.
 Also,
 $$ \pi'(\Gamma_2') \cap Z(H_p')
 \supset [\Gamma_2',\Gamma_2']
 = [\Gamma_2,\Gamma_2] ,$$ 
 so $\pi'(\Gamma_2') \cap Z(H_p')$ is a lattice in $[H_p,H_p] =
Z(H_p')$. Thus, we conclude that 
 $\pi'(\Gamma_2')$ is a lattice in~$H_p'$. Because $\pi'(\Gamma_2')$ is
contained in the arithmetic lattice $\pi'(\Gamma_2)$, this implies that
$\pi'(\Gamma_2')$ is arithmetic.

From the preceding paragraph, we know that $\mu_p$ is an isomorphism of
arithmetic lattices in $H'_p$. Let $\Fr \colon H \rightarrow H'_p$
denote the group isomorphism induced by the Frobenius automorphism $x
\rightarrow x^p$ of the ground field $\field $. Then there exist
arithmetic lattices $\hat{\Gamma _1}, \hat{\Gamma _2}$ in $H$, such
that $\Fr (\hat{\Gamma _1})=\Gamma '_1$ and $\Fr (\hat{\Gamma _2}) =
\pi '(\Gamma '_2)$, and an isomorphism
 $\lambda = \Fr ^{-1} \circ \mu _p \circ \Fr \colon \hat{\Gamma _1}
\rightarrow \hat{\Gamma _2}$.
 By Corollary~\ref{cor:heis}, we can virtually extend $\lambda $ to a
virtual automorphism $\Lambda $ of $H$. Then $\Lambda '_p= \Fr \circ
\Lambda \circ \Fr ^{-1}$ is a virtual automorphism of $H'_p$ virtually
extending $\mu _p$.

Let $\tilde{\Lambda }_p = \Lambda '_p \times \Psi$, so
$\tilde{\Lambda}_p$ is a virtual automorphism of $H_p$. We can define a
map $\zeta$ on some finite index subgroup of $\Gamma _1$ by $\zeta
(\gamma)=\lambda _p(\gamma )\tilde{\Lambda }_p(\gamma )^{-1}$. By
Lemma~\ref{lemma:ext}, $\zeta $ virtually extends to $\hat{\zeta}
\colon H_p \rightarrow Z(H_p)$. Then $\Lambda _p=\tilde{\Lambda }_p
\zeta$ is a virtual endomorphism of $H_p$. Since $\ker (\zeta )
\supset [H_p,H_p]$ we conclude (much as in the proof of
Corollary~\ref{G2rigid}) that $\Lambda _p$ is a virtual automorphism.
It is easy to see that it virtually extends $\lambda _p$.
 \end{proof}

\end{document}